\documentclass[12pt,reqno]{article}
\usepackage{amsfonts}
\usepackage[utf8x]{inputenc}
\usepackage[english]{babel}
\usepackage{amsmath,amssymb,amsfonts,wasysym}   
\usepackage{graphics}                           
\usepackage{amsmath}  
\usepackage{mathtools} 
\usepackage{hyperref}             
\usepackage{comment}
\newtheorem{theorem}{Theorem}[section]

\newtheorem{proposition}[theorem]{Proposition}

\newtheorem{observation}[theorem]{Observation}

\newcommand*{\Op}{\mathop{}\!\mathrm{Op}}
\newcommand*{\Proof}{\mathop{}\!\mathit{Proof}}

\usepackage{thmtools}
\usepackage{thm-restate}
\usepackage{hyperref}

\usepackage{cleveref}

\begin{document}
\title{On the Weyl transform and the quantization of the hypersphere}
\author{Simone Camosso$^*$}
\date{}
{\renewcommand{\thefootnote}{\fnsymbol{footnote}}
\setcounter{footnote}{1}
\footnotetext{\textbf{e-mail}: r.camosso@alice.it}
\setcounter{footnote}{0}
}

\maketitle

\begin{center}
\textbf{Abstract}
\end{center}
\noindent
After having dealt with the classical Weyl quantization, the deformation quantization and the recently (but old) Born-Jordan quantization, the purpose of the article is a sort of ``monomial quantization'' of $S^2$. The result of the impossibility of a rigorous quantization of the sphere is well known and treated in the literature, despite everything the case of the hydrogen atom remains one of the most interesting cases in the modeling of quantum theories.

\vspace{1cm}


\tableofcontents
\section{Introduction}
Let us consider $x$ and $p$ to be respectively the position and the momentum in a classical space.
It is well know the ``famous'' result due to Groenewold (\cite{Gr}) and van Hove (\cite{H1}) that showed that the general Dirac correspondence is not valid for all monomials. In fact there is a counterexample in dimension one for the operator $q^2p^2$. It is possible, through some calculations, discover two formulas for the quantization of the classical observable $x^2p^2$:

\begin{equation}
	\label{first_gr_formula11}
	\Op\left(x^{2}p^{2}\right)\,=\, \widehat{x}^{\,2}\widehat{p}^{\,2}-2i\hbar\widehat{x}\widehat{p}-\frac{2}{3}\hbar^2,
\end{equation}
\noindent 
and 
\begin{equation}
	\label{first_gr_formula22}
	\Op\left(x^{2}p^{2}\right)\,=\, \widehat{x}^{\,2}\widehat{p}^{\,2}-2i\hbar\widehat{x}\widehat{p}-\frac{1}{3}\hbar^2,
\end{equation}
\noindent 
that are in conflict. 

From this impossibility theorem we can try different rules in order to have a ``monomial quantization'' of the angular momentum on the sphere. In what follows we analyze the famous Weyl quantization introduced by Weyl with the beautiful property that bring the definition of moyal product. 
In a second part of the paper we analize another kind of quantization called the Born--Jordan quantization. 

At the end we try to find a ``monomial Born-Jordan quantization of the angular momentum on the hypersphere".

\section{The Dirac function and the Fourier transform}

The Dirac delta function is formally defined as:

\begin{equation}
\label{delta_1}
\delta(x)\,=\, \begin{cases}   +\infty \ \ \text{if} \ \ x=0 \\  \ \ \ \ 0 \ \ \text{if} \ \ x\not=0 \end{cases},
\end{equation}

with the property that:

\begin{equation}
	\label{delta_2}
	\int_{-\infty}^{+\infty}\delta(x)dx\,=\, 1.
\end{equation}

Paul Dirac defined this function as a ``generalized function'' (a distribution) in 1930. The problem is that the definition is not compatible with the condition $(\ref{delta_2})$. For this reason Dirac suggested to define $\delta(x)\,=\, \lim_{\varepsilon \rightarrow 0^{+}}\delta(x,\varepsilon)$, where $\delta(x,\varepsilon)$ is a generic function such that:

\begin{equation}
	\label{delta_3}
	\lim_{\varepsilon \rightarrow 0^{+}}\delta(x,\varepsilon)\,=\, \begin{cases}   +\infty \ \ \text{if} \ \ x=0 \\  \ \ \ \ 0 \ \ \text{if} \ \ x\not=0 \end{cases},
\end{equation}
\noindent 
and

\begin{equation}
	\label{delta_4}
	\int_{-\infty}^{+\infty}\delta(x,\varepsilon)dx\,=\, 1.
\end{equation}

An example of these functions is the Gaussian:

\begin{equation}
	\label{delta_5}
	\delta(x,\varepsilon)\,=\, \frac{1}{\varepsilon\sqrt{\pi}}e^{-\frac{x^2}{\varepsilon^2}}.
\end{equation}

If $x$ denote a vector in $\mathbb{R}^{n}$ instead an element of $\mathbb{R}$, the Dirac function $\delta(x)$ can be defined without different modifications respect the usual definition on $\mathbb{R}$. In what follows we will denote in the same way  by $\delta$ the $1$--Dirac distribution or the $n$--Dirac distribution.  

Let us denote the Fourier transform of a function $f(x)\in L^{1}(\mathbb{R}^{n})$ as:

\begin{equation}
	\label{fourier_1}
	\mathcal{F}(f)(\xi)\,=\,\int_{\mathbb{R}^n}f(x)e^{-2\pi ix\cdot \xi}dx.
\end{equation}
\noindent 
and the inverse:

\begin{equation}
	\label{fourier_2}
	\mathcal{F}^{-1}(f)(x)\,=\,\int_{\mathbb{R}^n}f(\xi)e^{2\pi ix\cdot \xi}d\xi.
\end{equation}

The delta function can be expressed in the ``sense of a distribution'' as:

\begin{equation}
	\label{delta_6}
	\delta(x)\,=\,\int_{\mathbb{R}^n}e^{2\pi ix\cdot \xi}d\xi.
\end{equation}

From the relation $(\ref{delta_6})$ it is possible to prove the following identity.

\begin{proposition}
	\begin{equation}
		\label{delta_7}
		f(x)\,=\,\int_{\mathbb{R}^n}\delta(x-y)f(y)dy.
	\end{equation}
	
\end{proposition}

$\Proof$ By $(\ref{delta_6})$ we have that:

$$
\begin{multlined}
	\int_{\mathbb{R}^n}\delta(x-y)f(y)dy\,=\, \int_{\mathbb{R}^n}\int_{\mathbb{R}^n}e^{2\pi i(x-y)\cdot \xi}d\xi f(y)dy\,=\, \\
	\int_{\mathbb{R}^n}e^{2\pi ix\cdot \xi}\int_{\mathbb{R}^n}e^{-2\pi iy\cdot \xi} f(y)dyd\xi \,=\,   		\int_{\mathbb{R}^n}e^{2\pi ix\cdot \xi}\mathcal{F}(f)(\xi)d\xi\,=\, \\ \,=\, \mathcal{F}^{-1}\left(\mathcal{F}(f)\right)(x)\,=\, f(x).
\end{multlined}
$$
\hfill $\Box$
	
In the next section we will revisit the Weil Quantization. This kind of quantization is a remarkable way to associate ``physical observables'' as $q$ the position and $p$ the momentum to ``quantum operators'' respectively $Q$ and $P$. We remember here the rules of commutation between these operators:

\begin{equation}
	\label{commutations}
	[P_{j},P_{k}]\,=\, [Q_{j},Q_{k}]\,=\, 0, \ \  \ \ [P_{j},Q_{k}]\,=\, \frac{\hbar}{i},
\end{equation}
\noindent 
for $k,j\,=\, 1, \ldots , n$ components of the momentum operator $P$ and position operator $Q$.

\section{The Weyl quantization}

The idea of the Weyl quantization is to associate to a function $\sigma(\xi,x)$ of $\xi,x\in\mathbb{R}^{n}$ an operator $\sigma(D,X)$ from the Schwartz space $\mathcal{S}(\mathbb{R}^{n})$ to the Schwartz distributional space $\mathcal{S}'(\mathbb{R}^{n})$. This association is given by the following Weyl operator:

\begin{equation}
	\label{Weyl_quant}
	\sigma(D,X)f(x)\,=\, \int_{\mathbb{R}^{n}} \int_{\mathbb{R}^{n}}  \sigma\left(\xi, \frac{x+y}{2}\right)e^{2\pi i(x-y)\xi}f(y)dyd\xi,
\end{equation}
\noindent 
for some $f(x)\in \mathcal{S}(\mathbb{R}^{n})$. 

The integral kernel of $(\ref{Weyl_quant})$ si given by:

\begin{equation}
	\label{Weyl_quant_1.5}
 K_{\sigma}(x,y)\,=\, \int_{\mathbb{R}^{n}} \sigma\left(\xi, \frac{x+y}{2}\right)e^{2\pi i(x-y)\xi}d\xi,
\end{equation}
\noindent 
for this reason the formula $(\ref{Weyl_quant})$ can be written also as:

\begin{equation}
	\label{Weyl_quant_2}
	\sigma(D,X)f(x)\,=\, \int_{\mathbb{R}^{n}} K_{\sigma}(x,y)f(y)dy.
\end{equation}

In order to prove the formula $(\ref{Weyl_quant})$ the first step is to consider the inverse Fourier transform of the Fourier transform of $\sigma$:

\begin{equation}
	\label{proof-1}	\sigma(D,X)f(x)\,=\, \int_{\mathbb{R}^{n}} \int_{\mathbb{R}^{n}}  \mathcal{F}(\sigma)\left(p,q\right)\left(e^{2\pi i(pD+qX)}f\right)(x)dpdq.
\end{equation}

We can evaluate the quantity $e^{2\pi i(pD+qX)}f(x)$ using the following identity:

\begin{equation}
	\label{quantity}
	e^{2\pi i(pD+qX)}f(x)\,=\, e^{2\pi i qx + i\pi pq}f(x+p).
\end{equation}

We find that:

\begin{equation}
	\label{proof-2}
\sigma(D,X)f(x)\,=\, \int_{\mathbb{R}^{n}} \int_{\mathbb{R}^{n}} \int_{\mathbb{R}^{n}} \int_{\mathbb{R}^{n}}  \sigma\left(\xi,w\right)e^{-2\pi i(p\xi+qw)} e^{2\pi i qx + i\pi pq}f(x+p)d\xi dwdpdq.
\end{equation}

We observe that $\delta\left(x-w+\frac{p}{2}\right)\,=\, \int_{\mathbb{R}^{n}}e^{2\pi i \left(x +\frac{p}{2}-w\right)q}$ and $(\ref{proof-2})$ can be written as:

\begin{equation}
	\label{proof-3}
\sigma(D,X)f(x)\,=\, \int_{\mathbb{R}^{n}} \int_{\mathbb{R}^{n}} \int_{\mathbb{R}^{n}}  \sigma\left(\xi,w\right)e^{-2\pi ip\xi} \delta\left(x-w+\frac{p}{2}\right)f(x+p)d\xi dwdp.
\end{equation}

Now we observe that $\sigma\left(\xi,x+\frac{p}{2}\right)\,=\,\int_{\mathbb{R}^{n}}\sigma(\xi,w)\delta\left(x-w+\frac{p}{2}\right)dw$ and, substituting in $(\ref{proof-3})$ we deduce that:

\begin{equation}
	\label{proof-4}
\sigma(D,X)f(x)\,=\, \int_{\mathbb{R}^{n}} \int_{\mathbb{R}^{n}}\sigma\left(\xi,x+\frac{p}{2}\right) e^{-2\pi ip\xi} f(x+p)d\xi dp.
\end{equation}

If we set $x+p=y$, $dp=dy$ and $(\ref{proof-4})$ is transformed to $(\ref{Weyl_quant})$:

\begin{equation}
	\label{Weyl_quant_lasted}
	\sigma(D,X)f(x)\,=\, \int_{\mathbb{R}^{n}} \int_{\mathbb{R}^{n}}  \sigma\left(\xi, \frac{x+y}{2}\right)e^{2\pi i(x-y)\xi}f(y)dyd\xi.
\end{equation}

\begin{observation}
Details of the proof can be found in \cite{uno}, in particular the equality $(\ref{quantity})$ follows from the PDE problem:

\begin{equation}
	\label{PDE}
	\begin{cases} \partial_{t}g(x,t) =2\pi i (pD+qX)g(x,t) \\ g(x,0)\,=\, f(x) \end{cases},
\end{equation}
\noindent 
where $g(x,t)\,=\, \left(e^{2\pi i t(pD+qX)}f\right)(x)$. 
It is possible to reduce the problem solving an ODE observing that $ \partial_{t}g(x,t) -2\pi ipDg(x,t)=2\pi i qxg(x,t)$ is the directional derivative along the vector $(-p,1)$. Let us to set $x(t)=x-tp$ and $G(t)=g(x(t),t)$, we have the ODE problem:

\begin{equation}
	\label{ODE}
	\begin{cases} G'(t)\,=\, 2\pi i q(x-tp)G(t) \\ G(0)\,=\, f(x) \end{cases}.
\end{equation}

The problem $(\ref{ODE})$ has a solution: $g(x-tp,t)\,=\,G(t)\,=\, e^{2\pi i tqx - \pi it^2pq}$. Now if we set $t=1$ and $x\mapsto x+p$ we find the identity $(\ref{quantity})$. 

\end{observation}

Let us consider $\sigma, \tau \in\mathcal{S}(\mathbb{R}^{n})$, it is possible to consider the following twisted product:

\begin{equation}
	\label{product1}
	\sigma \# \tau (D,X) \,=\, \sigma(D,X)\tau(D,X).
\end{equation}

Let us consider in detail the previous formula:

\begin{equation}
	\label{product2}
	\begin{multlined}
	\sigma \# \tau (D,X) f(x) \,=\, \int_{\mathbb{R}^{4n}} \sigma\left(\zeta, \frac{x+z}{2}\right)\tau\left(\eta, \frac{y+z}{2}\right)e^{2\pi i(x-z)\zeta+(z-y)\eta}f(y)dyd\eta dz d\zeta \,=\, \\ \,=\, \int_{\mathbb{R}^{n}}K(x,y)f(y)dy.
	\end{multlined}
\end{equation}

By the Wigner transform:

\begin{equation}
	\label{product3}
	\begin{multlined}
		\sigma \# \tau (\xi,x) \,=\, \int_{\mathbb{R}^{n}} K\left(x+\frac{p}{2},x-\frac{p}{2}\right)e^{-\pi i \xi p} dp \,=\,  \\ \,=\, \int_{\mathbb{R}^{4n}}\sigma\left(\zeta, \frac{x+z}{2}+\frac{p}{4}\right)\tau\left(\eta, \frac{x+z}{2}-\frac{p}{4}\right)e^{2\pi i (x\zeta-z\zeta+z\eta-x\eta-\xi p)+\pi i p (\zeta+\eta)}d\eta dz d\zeta dp.
	\end{multlined}
\end{equation}

If we set $u= \frac{x+z}{2}+\frac{p}{4}$ and $v=\frac{x+z}{2}-\frac{p}{4}$, we have that $dz dp =4^n dudv$ and $x\zeta-z\zeta+z\eta-x\eta-\xi p+ \frac{1}{2}p (\zeta+\eta)=2\omega\left((\xi-\eta,x-v),(\xi-\zeta,x-u)\right)$, where $\omega$ is the symplectic form on $\mathbb{R}^{2n}$. Thus we have that:

\begin{equation}
	\label{product4}
	\begin{multlined}
		\sigma \# \tau (\xi,x) \,=\, 4^n\int_{\mathbb{R}^{4n}}\sigma(\zeta,u)\tau(\eta,v)e^{4\pi i \omega\left((\xi-\eta,x-v),(\xi-\zeta,x-u)\right)}d\zeta d u d\eta dv.
	\end{multlined}
\end{equation}

Now by the Taylor formula we have that:

\begin{equation}
	\label{product5}
	\begin{multlined}
		\sigma(\zeta, u)\tau(\eta,v) \,=\, \sum_{|\alpha|+|\beta|\leq k}\frac{\partial_{x}^{\alpha}\sigma(\zeta,x)\partial_{x}^{\beta}\tau(\eta,x)}{\alpha!\beta!}(u-x)^{\alpha}(v-x)^{\beta} + R
	\end{multlined}
\end{equation}
\noindent 
where $R$ is a remainder. We observe that:

$$(u-x)^{\alpha}(v-x)^{\beta}e^{4\pi i \left[(x-u)(\xi-\eta)-(x-v)(\xi-\zeta)\right]}\,=\, \frac{(-1)^{|\beta|}}{2^{|\alpha|+|\beta|}}D_{\eta}^{\alpha}D_{\zeta}^{\beta}e^{4\pi i \left[(x-u)(\xi-\eta)-(x-v)(\xi-\zeta)\right]}.$$

Now plugging it in the formula $(\ref{product4})$, after integration by parts and the application of the Fourier inversion formula, we find that:

\begin{equation}
	\label{product6}
	\begin{multlined}
	\sigma \# \tau (\xi,x)\,=\,\sum_{|\alpha|+|\beta|\leq k}\frac{(\pi i)^{|\alpha|+|\beta|}(-1)^{|\alpha|}}{\alpha!\beta!}D_{\xi}^{\beta} D_{x}^{\alpha}\sigma(\xi,x)D_{\xi}^{\alpha} D_{x}^{\beta}\tau(\xi,x) +R_{k}(\xi,x)\,=\, \\ \,=\, \sum_{j=0}^{k}\frac{(\pi i)^j}{j!}\left(D_{\xi,\sigma}D_{x,\tau}-D_{\xi,\tau}D_{x,\sigma}\right)^{j}\sigma(\xi,x)\tau(\xi,x).
	\end{multlined}
\end{equation}
\noindent 
where for example the notation $D_{\xi,\sigma}$ is a differentiation applied only to $\sigma$.

Now defining the momentum operator as $P=hD$, where $h$ is the Planck constant, we have that:

\begin{equation}
	\label{product7}
	\begin{multlined}
		\sigma \#_{h} \tau (\xi,x)\,=\, \sum_{j=0}^{k}\frac{(\pi ih)^j}{j!}\left(D_{\xi,\sigma}D_{x,\tau}-D_{\xi,\tau}D_{x,\sigma}\right)^{j}\sigma(\xi,x)\tau(\xi,x).
	\end{multlined}
\end{equation}
\noindent 
and we have that:

\begin{equation}
	\label{product8}
	\begin{multlined}
	\{\sigma, \tau\}	\,=\, \frac{2\pi i }{h}\left(\sigma \#_{h} \tau - \tau  \#_{h} \sigma\right)+\mathcal{O}(h^2),
	\end{multlined}
\end{equation}
\noindent 
where $\{\cdot,\cdot\}$ are the Poisson backet and the product $\#_{h}$ is also called the Moyal product. The formula $(\ref{product8})$ can be summarize by the following formula:

\begin{equation}
	\label{product9}
	\begin{multlined}
		\{\sigma, \tau\}	\,=\, \frac{2\pi i }{h}[\sigma,\tau]_{\#_{h}}+\mathcal{O}(h^2),
	\end{multlined}
\end{equation}
\noindent
that represents the Weyl correspondence as $h\rightarrow 0$. The correspondence will say that the Dirac quantum condition is satisfied asymptotically when $\hbar \rightarrow 0$. 

\section{The Weyl quantization on manifolds}

The Weyl quantization can be generalized on a manifold $M$. Let us consider $s$ a symbol in $S^{m}(M)$, then the Weyl quantization is given by the following formula:

\begin{equation}
	\label{Weyl_quant_manifold}
	[W(s)f](x)\,=\, \frac{1}{2\pi \hbar}\int_{T^{\vee}M}\chi(x,y) e^{\frac{i}{\hbar}g(\exp^{-1}_{y}(x),\xi)}s\left(\tau_{g_{\frac{1}{2}}(x,y),y}\xi\right)f(y)dyd\xi,
\end{equation}
\noindent 
where $\chi$ is a properly supported cut--off function around the diagonal $M\times M$ such that $\exp^{-1}_{y}(x)$ is well defined for all $(x,y)$ in the support. The quantity $g_{\frac{1}{2}}(x,y)$ is the half  midpoint between $x$ and $y$ or $\exp{\frac{1}{2}(\exp^{-1}_{x}(y))}$. The notation $\tau_{x,y}$ denotes the parallel transport in $T^{\vee}M$ from $T^{\vee}_{y}M$ to $T^{\vee}_{x}M$ along the geodesic joining $x$ and $y$. In the end $g(\cdot, \cdot)$ is the metric on $M$. For detalis see \cite{due}.

The remarkable fact is that the formula $(\ref{Weyl_quant_manifold})$ gives the Weyl correspondence:

\begin{equation}
	\label{product9}
	\begin{multlined}
	[W(s),W(t)]	\,=\,-i\hbar W\left(\{s, t\}\right) +\mathcal{O}(h^2),
	\end{multlined}
\end{equation}
\noindent 
for every symbol $s,t\in S^{m}(M)$. Moreover the Weyl quantization induces a deformation quantization of the cotangent bundle $TM$. 

\section{The Weyl quantization of $l$ and the angular ``momentum dilemma''}

Let us consider the classical angular momentum operator of the 1s state. It is given by:

\begin{equation}
	\label{l_operator1}
	l\,=\, \left(x_{2}p_{3}-x_{3}p_{2},x_{3}p_{1}-x_{1}p_{3},x_{1}p_{2}-x_{2}p_{1}\right),
\end{equation}
\noindent 
where $x \,=\, \left(x_{1},x_{2},x_{3}\right)\in\mathbb{R}^{3}$ and $p \,=\, \left(p_{1},p_{2},p_{3}\right)\in\mathbb{R}^{3}$. Let us consider only the third component of $l$, denoted by $l_{3}\,=\, x_{1}p_{2}-x_{2}p_{1}$. Let us consider the operator $\widehat{l}_{3}\,=\, \widehat{x}_{1}\widehat{p}_{2}-\widehat{x}_{2}\widehat{p}_{1}$ and the square of this component:

\begin{equation}
	\label{squareroot_l_squared}
	\widehat{l}_{3}^{\,\,2}\,=\, \widehat{x}_{1}^{\,2}\widehat{p}_{2}^{\,2}-\widehat{x}_{2}^{\,2}\widehat{p}_{1}^{\,2}-\widehat{x}_{1}\widehat{p}_{1}\widehat{p}_{2}\widehat{x}_{2}-\widehat{p}_{1}\widehat{x}_{1}\widehat{x}_{2}\widehat{p}_{2}.
\end{equation}

Let us consider the Weyl quantization rule for monomials $x^{m}p^{n}$:

\begin{equation}
	\label{Weyl_quantization_rule}
	\Op_{\text{Weyl}}\left(x^{m}p^{n}\right)\,=\, \frac{1}{2^{m}}\sum_{k=0}^{m}\binom{m}{k}\widehat{x}^{\,k}\widehat{p}^{\,n}\widehat{x}^{\,m-k}\,=\, \frac{1}{2^{n}}\sum_{k=0}^{n}\binom{n}{k}\widehat{p}^{\,k}\widehat{x}^{\,m}\widehat{p}^{\,n-k},
\end{equation}
\noindent 
and we apply the rule ($\ref{Weyl_quantization_rule}$) in order to quantize $l^{2}$ (it is possible to quantize only $l_{3}$ because for $l_{1}$ and $l_{2}$ the procedure it is analogue). What happen is that the difference between the Weil transform of $l^2$ and $\widehat{l}^{\,\,2}$ is $\frac{3}{2}\hbar^2$. This discrepacy is called the ``momentum dilemma'' and a deep discussion can be found in \cite{DS} and \cite{G1}.

\section{The Shubin $\tau$--ordering for monomials}
Let $n,m$ be two non--negative integers. Let us consider the following two rules:

\begin{equation}
	\label{one_S}
	x^{m}p^{n}\mapsto \sum_{k=0}^{m}\binom{m}{k}\tau^k\left(1-\tau\right)^{m-k}\widehat{x}^{\,k}\widehat{p}^{\,n}\widehat{x}^{\,m-k}
\end{equation}
\noindent 
or equivalently 

\begin{equation}
	\label{one_S2}
	x^{m}p^{n}\mapsto \sum_{k=0}^{n}\binom{n}{k}\tau^{n-k}\left(1-\tau\right)^{k}\widehat{p}^{\,k}\widehat{x}^{\,m}\widehat{p}^{\,n-k},
\end{equation}
\noindent 
where $\widehat{x}^{\,j}$ and $\widehat{p}^{\,j}$ are the usual quantum mechanic operators that on a wave function $\psi$ are given by: $\widehat{x}^{\,j}\psi(x)\,=\, x^j\cdot \psi(x)$ and $\widehat{p}^{\,j}\psi(x)\,=\, \left(-i\hbar\partial_{x}\right)^{\,j} \psi(x)$ for some $j$. In the setting of the Shubin $\tau$--ordering $\tau$ is a real parameter (see \cite{S}). The rules $(\ref{one_S})$ and $(\ref{one_S2})$ are a sort of generalization of the Weyl quantization ($\tau=\frac{1}{2}$).

\section{The Born--Jordan quantization}

The Born--Jordan quantization rule (BJ) is given by the integration of ($\ref{one_S2}$) on the interval $[0,1]$ respect to the real parameter $\tau$. The integral must be understood in some ``appropriate sense'' and the calculations involves the following integral:

\begin{equation}
	\label{imp_integral}
	\int_{0}^{1}\tau^{n-k}\left(1-\tau\right)^{k}d\,\tau\,=\, \frac{(n-k)!k!}{(n+1)!}.
\end{equation}

Let us denote the BJ operator as:

\begin{equation}
	\label{imp_integral}
	\Op_{\text{BJ}}\left(x^{m}p^{n}\right)\,=\, \int_{0}^{1}\Op_{\text{Shubin}}\left(x^{m}p^{n}\right)d\,\tau.
\end{equation}

The major reference is the book of \cite{G2}. The integral gives the following rule:
\begin{equation}
	\label{BJ_quantization_rule}
	\Op_{\text{BJ}}\left(x^{m}p^{n}\right)\,=\, \frac{1}{m+1}\sum_{k=0}^{m}\widehat{x}^{\,m-k}\widehat{p}^{\,n}\widehat{x}^{\,k}.
\end{equation}

\section{Difference between Born--Jordan and Weyl quantization}

Let us denote with BJ  the Born Jordan rule $(\ref{BJ_quantization_rule})$ and with W the Weyl rule $(\ref{Weyl_quantization_rule})$. Let us consider the monomial $x^{m}p^{n}$ ($n,m \in \mathbb{N}$ different from $0$), then the two rules coincide if $m+n\leq 2$. In general they are different, for example if $n=m=2$ we can prove the following proposition.

\begin{proposition}
The BJ quantization and the W quantization of the monomial $x^2p^2$ differ by the quantity $\frac{\hbar^2}{2}$.
\end{proposition}	

$\Proof .$ It is a simple calculation using several times the commutation relation $\left[\widehat{x},\widehat{p}\right]\,=\, i\hbar$.

\hfill $\Box$

\begin{proposition}
	The BJ quantization of the angular momentum operator $l^2$ for the 1s orbit is:
	
\begin{equation}
\label{bj_atom_1s}
\Op_{\text{BJ}}(l^2)\,=\, \widehat{l}^{\,\,2}+2\hbar^2.
\end{equation}
\end{proposition}	

$\Proof .$ 
Let us denote with $l^2\,=\, l_{1}^{2}+l_{2}^{2}+l_{3}^{2}$ the length of the angular momentum vector and, the third component $ l_{3}$. 
Let us consider the symmetrization rule ($\ref{one_S}$) for the $l_{3}$ classical function and the definition of the Born--Jordan operator ($\ref{imp_integral}$), we have that:

$$\Op_{\text{BJ}}(l^2_{3})\,=\, \widehat{x}_{1}^{\,2}\widehat{p}_{2}^{\,2}+\widehat{x}_{2}^{\,2}\widehat{p}_{1}^{\,2}-\frac{2}{3}\widehat{x}_{1}\widehat{p}_{1}\widehat{x}_{2}\widehat{p}_{2} - \frac{1}{3}\widehat{x}_{1}\widehat{p}_{1}\widehat{p}_{2}\widehat{x}_{2}-\frac{1}{3}\widehat{p}_{1}\widehat{x}_{1}\widehat{x}_{2}\widehat{p}_{2}-\frac{2}{3}\widehat{p}_{1}\widehat{x}_{1}\widehat{p}_{2}\widehat{x}_{2}.$$

Now, using the relations $\left[\widehat{x}_{k},\widehat{p}_{k}\right]$ (for $k=1,2$) we find that:

$$\Op_{\text{BJ}}(l^2_{3})\,=\,\widehat{l}_{3}^{\,\, 2}+\frac{2}{3}\hbar^{2}.$$

The result follows considering also the other components of $l^2$.

\hfill $\Box$

In general Born quantization and Weyl quantization differ by $\frac{\hbar^2}{2}$ factor.

\section{A monomial deformed Born--Jordan quantization of $l$ for $S^3$}

In analogy with the case of a $3$--dimensional space we can define the angular momentum for $x,p\in\mathbb{R}^{4}$. In this case we have the matrix:

\begin{equation}
	\label{quadri_angular_momentum}
	l\,=\, \left(\begin{array}{cccc} 0 & l_{12} & l_{13} & l_{14} \\ -l_{12} & 0 & l_{23} & l_{24} \\ -l_{13} & -l_{23} & 0 & l_{34} \\ -l_{14} & -l_{24} &-l_{34} & 0 \end{array}\right)
\end{equation}
\noindent 
where $l_{ij}\,=\,x_{i}p_{j}-x_{j}p_{i}$ for $i,j\,=\,1,2,3,4$.  

We can define the operator $l^2\,=\, \frac{1}{2}\sum_{i,j=1}^{4}l_{ij}^2$. The Weyl quantization brings to have a quantization value of $3\hbar^2$. Generalizing the previous resut we have that The Born Jordan quantization brings a value of $4\hbar^2$. The difference between the two methods is $\hbar^2$. We did a sort of ``monomial quantization'' of the hypersphere. Conjecture: for a sphere of dimension $n$ is the Born Jordan quantization be $2(n-2)\hbar^2$ ?

\end{document}